\documentclass[12pt]{article}
\usepackage{amsfonts}
\usepackage{amsfonts}
\usepackage{amsfonts}
\usepackage{amsfonts}
\usepackage{pdflscape}
\usepackage[hidelinks]{hyperref}

\usepackage{amsthm}
\usepackage{amsmath,amssymb}

\usepackage{adjustbox}
\usepackage{graphicx}

\allowdisplaybreaks[4]

\numberwithin{equation}{section}

\newtheorem{definition}{Definition}

\title{A Note on the Control of Mortality for \\Nested Recurrence Relation $H(n)=H(n-H(n-2))+H(n-H(n-3))$}
{\footnotesize
\author{
{\sc Altug Alkan}$^a$\footnote{All correspondences should be addressed to: \href{mailto:altug.alkan@pru.edu.tr}{\texttt{altug.alkan@pru.edu.tr}}},
{\sc Orhan Ozgur Aybar}$^b$ and {\sc Zehra Akdeniz}$^c$
\\
\\
$^a${\small Piri Reis University, Institute of Graduate Studies}\\
$^b${\small Piri Reis University, Continuing Education Center} \\
$^c${\small Piri Reis University, Faculty of Science and Letters} \\
{\small 34940 Tuzla/Istanbul, Turkey}
\\
\,{}\\
{}
}

\date{\today}

\begin{document}
\maketitle

\begin{abstract}
\vskip 0.1cm

In this work, we investigate $H$-recurrence which is defined by $H(n) = H(n-H(n-2)) + H(n-H(n-3))$ thanks to a recent approach to certain recurrences such as Conway and Hofstadter $Q$ recursions.


\vskip 0.1cm
\end{abstract}

\baselineskip 15pt

\parskip 5pt

\section{Introduction and Motivation}
It is known that meta-Fibonacci recurrence relations are considered one of the most resistant subclasses of nonlinear recurrences in terms of common analysis and proof techniques of difference equations~\cite{1}. While there is no universal agreement on generalized definition of nested recurrence
relations, concept of meta-fibonacci has been introduced by Douglas Hofstadter with the invention of original $Q$-sequence (A005185 in OEIS~\cite{2}) is defined by $Q(n) = Q(n-Q(n-1)) + Q(n-Q(n-2))$ with initial conditions $Q(1)=Q(2)=1$~\cite{3}, see Figure 1 for its curious scatterplot. Although very little is known about mathematical properties of $Q$-sequence, many related recurrences and solution sequences have been examined and the literature on this research field has rigorous results especially for slow and quasi-periodic solutions [4-10]. At the same time, there are curios examples that belong to various
experimental mathematics studies provide meaningful insigths for not only perfect solutions but also unpredictable sequences that are generated by meta-Fibonacci recurrences [11-14]. At this direction, we focus on $H(n) = H(n-H(n-2)) + H(n-H(n-3))$ recurrence that is investigated recently and meaningful results are proved in terms of the existence of new kind of solution families which are related to systems of nested recurrences that resemble Golomb's recurrence $G(n)=G(n-G(n-1))+1$~\cite{14}. As we try to study $H(n)$ for a future research direction of previous paper~\cite{14}, we will adopt the approach of generalizations introduced earlier~\cite{15,16,17}. Some of examples with definitions as below see Figure 2, 3, 4 for corresponding generalizations.

%indeed in mathematical In this work, we connect some recent results that computational perspective on initial condition patterns that are formulated by asymptotic properties provide a variety of mathematical discoveries on some famous nested recurrence relations. We mainly focus on certain sequences that are related by Conway and Hofstadter-like recurrences for motivation.

%\begin{figure}[!htbp]
%\minipage{0.52\textwidth}
  %\includegraphics[width=\linewidth]{Hofterms1.png}
%\endminipage\hfill
%\minipage{0.5\textwidth}
  %\includegraphics[width=\linewidth]{qbro.png}
%\endminipage\hfill
%\caption {Scatterplots of Hofstadter's $Q$-sequence and $Q_{b}(n)$}
%\label{fig:Q}
%\end{figure}

\begin{figure}[!htbp]
\begin{center}
\includegraphics[width=0.58\textwidth]{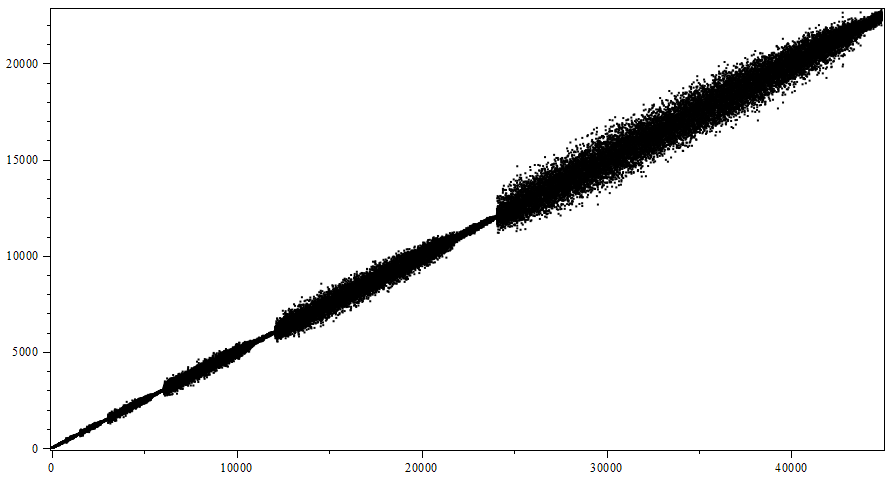}
\caption{Scatterplot of Hofstadter's $Q$-sequence for initial generations.}
\label{fig:Q}   
\end{center}
\end{figure}

%This paper is structured as follows.  In Section~\ref{sec:gen}, we review that the existing generalization of Hofstadter's $Q$-sequence according to the initial condition formulation. Then, in Section~\ref{sec:2analyse} and Section~\ref{sec:3analyse}, selected members of this curious sequence family are mentioned in detail based on their generational structures with the computational perspective.  In Section 4, we remember the limits of such generalization attempt while in Section 5 we observe the results of same approach to well-known slow families of solutions. Finally, some concluding remarks are offered in Section~\ref{sec:conc}.

%\newpage

%\newpage
%%%%%%%%%%%%%%%%%%%%%%%%%%%%%%%%%%%%%%%%%%%%%%%%%%%%%%%%%%%%%%%%%%%%%%%%%%%

%\section{Some Generalizations On Certain Recurrences}\label{sec:gen}
Definition 1 is related to generalization of original $Q$-sequence. See Figure 2 in order to observe examples $Q_{4,2}(n)$ and $Q_{5,2}(n)$ sequences which are similar conjectural properties with original $Q$-sequence~\cite{15}.

\begin{definition}

Let $Q_{d,\ell}(n)$ be defined by the recurrence $Q_{d,\ell}(n) = \sum\limits_{i=1}^\ell Q_{d,\ell}(n$ $-$ $Q_{d,\ell}(n$ $-$ $i))$ for $n>d \cdot \ell$, $\ell\geq 2$, $d\geq 1$, with the initial conditions $Q_{d,\ell}(n) = \lceil \frac{n\cdot(\ell-1)}{\ell}\rceil$ for $1\leq n\leq d\cdot\ell$.

\end{definition}

\begin{figure}[!htbp]
\minipage{0.5\textwidth}
  \includegraphics[width=\linewidth]{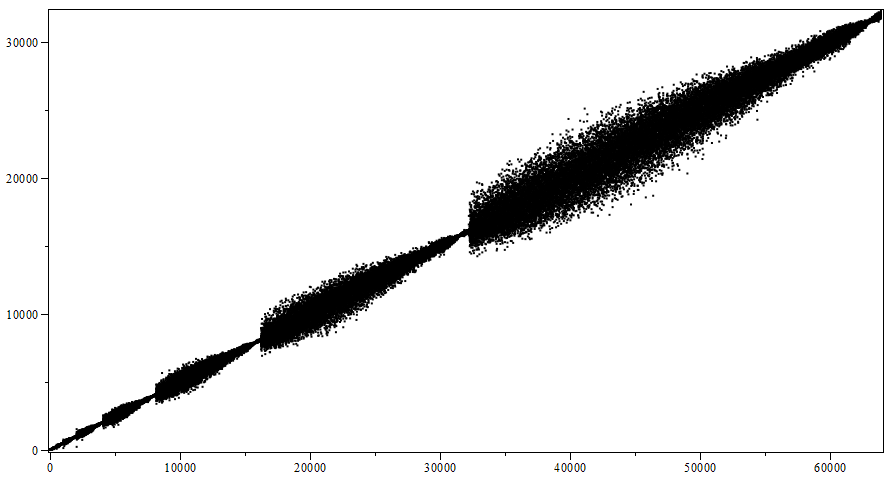}
\endminipage\hfill
\minipage{0.5\textwidth}
  \includegraphics[width=\linewidth]{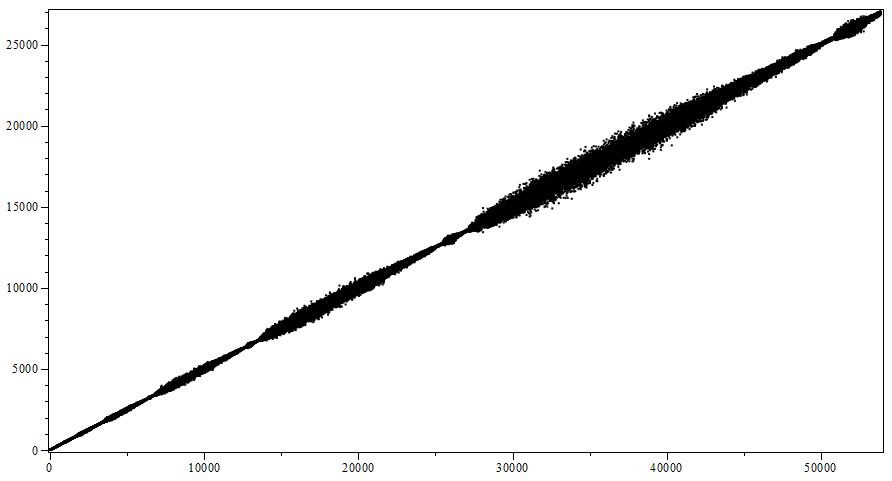}
\endminipage\hfill
\caption {Scatterplots of $Q_{4,2}(n)$ and $Q_{5,2}(n)$ for initial generations.}
\label{fig:Qx}
\end{figure}

\newpage

Definition 2 and 3 are related to generalization of Hofstadter-Conway $\$10000$ sequence, see Figure 3 and 4 in order to observe fractal-like collections of fractal-like sequences based on curious generational structures ~\cite{16,17}. These results which generate highly similar behaviour on solution sequences provide motivation to conduct experiments on similar intial condition patterns for $H(n)$ recurrence because  if $\lim_{n\to\infty}\frac{H(n)}{n}$ exists, it must be equal to $\frac{1}{2}$.
\begin{definition}

Let $c_{i}(n) = c_{i}(c_{i}(n-1)) + c_{i}(n-c_{i}(n-1))$ for $n > 4\cdot i $, with the initial conditions $c_{i}(n) = \lceil \frac{n}{2}\rceil$ for $1 \leq n \leq 4\cdot i$.

\end{definition}

\begin{figure}[!htbp]
\begin{center}
\includegraphics[width=0.8\textwidth]{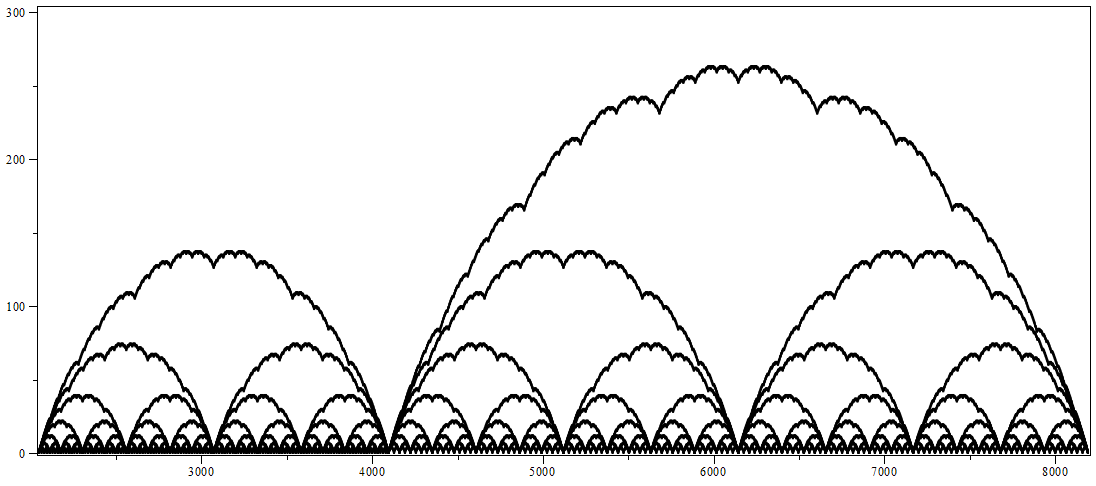}
\caption{Scatterplot of $c_{2^t}(n) - \frac{n}{2}$  for $2^{11} \leq n \leq 2^{13}$ and $0 \leq t \leq 7$.}
\label{fig:Q}   
\end{center}
\end{figure}
%newpage
%\begin{definition}

%Let $c^{*}_{i}(n) = n-c^{*}_{i}(c^{*}_{i}(n-1)) - c^{*}_{i}(n-c^{*}_{i}(n-1))$ for $n > 2^i$, with the initial conditions $c^{*}_{i}(n) = \lceil \frac{n}{2}\rceil$ for $1 \leq n \leq 2^{i}$.

%\end{definition}

%\begin{figure}[!htbp]
%\begin{center}
%\%includegraphics[width=0.99\textwidth]{Con2.png}
%\caption{Scatterplot of $c^{*}_{i}(n)- \frac{n}{2}$ for initial members and generations}
%\label{fig:Q}   
%\end{center}
%\end{figure}

%\begin{definition}

%Let $b_{i}(n) = n-b_{i}(b_{i}(n-i)) - b_{i}(n-b_{i}(n-i))$ for $n > 2 \lceil \frac{i}{2}\rceil$, with the initial conditions $b_{i}(n) = \lceil \frac{n}{2}\rceil$ for $1 \leq n \leq 2 \lceil \frac{i}{2}\rceil$.

%\end{definition}

\begin{definition}

Let $a_{i,j}(n) = a_{i,j}(a_{i,j}(n-j)) + a_{i,j}(n-a_{i,j}(n-1))$ for $n > 2 (i-1) + 3  \lfloor \phi  i \rfloor $, with the initial conditions $a_{i,j}(n) =\left \lfloor \frac{n+2}{1+\phi} \right \rfloor$ for $1 \leq n \leq  2 (i-1) + 3  \lfloor \phi  i \rfloor $ where $\phi = \frac{1+\sqrt{5}}{2}$ and $j \in \{1,2\}$. 

\end{definition}

\begin{figure}[!htbp]
\begin{center}
\includegraphics[width=0.8\textwidth]{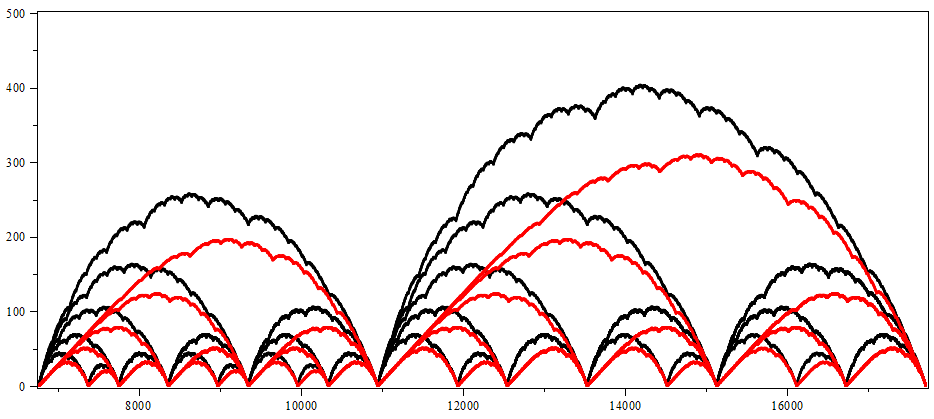}
\caption{Scatterplots of $a_{i,1}(n)- \left \lfloor \frac{n+2}{1+\phi} \right \rfloor$ (black) and $a_{i,2}(n)- \left \lfloor \frac{n+2}{1+\phi} \right \rfloor$ (red) for $F(20) \leq n \leq F(22)$ and $1 \leq i \leq 8$ where $F(n)$ is $n^{\text{th}}$ Fibonacci number.}
\label{fig:Q}   
\end{center}
\end{figure}

%\begin{figure}[!htbp]
%\begin{center}
%\includegraphics[width=0.99\textwidth]{conways.eps}
%\caption{Scatterplot of $Q_{2,2}(n)$}
%\label{fig:1122}   
%\end{center}
%\end{figure}

\newpage

\newpage
\section{An Analysis of Mortality for $H(n)$ Recurrence }
%\subsection{Selected Solutions to $Q_{d,2}(n)$}\label{sec:2analyse}
Because of the main theme mentioned in the previous section, we construct a generalization as below in order to search a long chaotic, perhaps infinite, sequence based on a generational structure, see Table 1 and Figure 5.

\begin{definition}

Let $L(k)$ is the number of terms of $H_{k}(n)$ sequence if $H_{k}(n)$ is finite, where $H_{k}(n) = H_{k}(n-H_{k}(n-2)) + H_{k}(n-H_{k}(n-3))$ for $n > k $, with the initial conditions $H_{k}(n) = \lceil \frac{n}{2}\rceil$ for $1 \leq n \leq k$. $L(k) = - 1$ if $H_{k}(n)$ is infinite sequence.

\end{definition}

%We heuristically expect its solutions of $H$-recurrence to fluctuate around the line $y=\frac{x}{2}$. Hence, we prefer to study the corresponding noise sequence $S_{h}(n) = H_{114}(n)-\frac{n}{2}$.
%Approximate self-similar block structures of certain members of $Q_{d,2}(n)$ family can be studied thanks to auxilary sequences similar with different works which give definitions of generations~\cite{1, 5,17,18}. 

%\begin{con}
%$c(n) - \frac{n}2$ $\ge$ $\left\lvert c^*(n) - \frac{n}2 \right\rvert$ for all $n \ge 1$
%\end{con}

\begin{table*}[!htbp]
\begin{center}

\begin{adjustbox}{max width=1\textwidth}
\begin{tabular}{cccccc}
\noalign{\smallskip}
$ $ & $ $ & $ $ & $m$ & $ $ & $ $\\
\noalign{\smallskip}\hline\noalign{\smallskip}
\noalign{\smallskip}
$ $ & $1$ & $2$ & $3$ & $4$ & $5$\\
\noalign{\smallskip}\hline\noalign{\smallskip}
$L(m+2)$ &53	&42	 &265	 &24 &39\\
$L(m+7)$ &1399	 &13270 &22308 &679&5665\\
$L(m+12)$ &7042 &860 &2214	&23110	&6682\\
$L(m+17)$ &52548&1539	&257156	&1001847	&49675\\
 $L(m+22)$ &69969&620771& 717281  &426570&86781\\
 $L(m+27)$ & 1846835&505915& 88391 &465589&488724\\
\hline
\noalign{\smallskip}
\end{tabular}
\end{adjustbox}
\caption {The values of $L(n)$ sequence for $n \leq 32$.}
\label{tab:table3}
\end{center}
\end{table*}

\begin{figure}[!htbp]
\begin{center}
\includegraphics[width=0.6\textwidth]{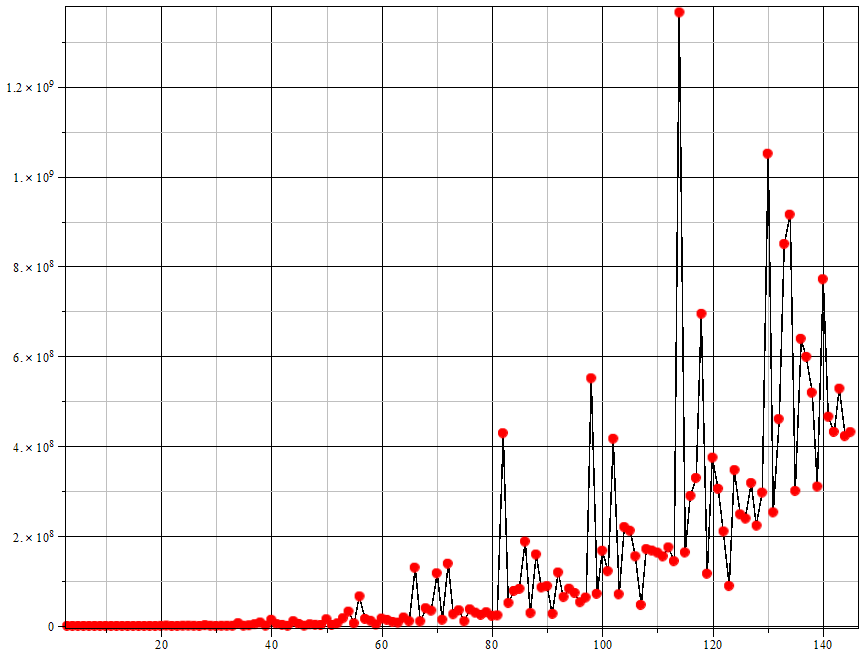}
\caption{Plot of $L(n)$ for $n \leq 145$. In particular, $L(114) = 1365203701$ is the maximum value in this range.}
\label{fig:Q}   
\end{center}
\end{figure}

Since $L(114) = 1365203701$ is the maximum value for $n \leq 145$ while $L(146) = 1747189782$, so we can try to search a generational structure in order to depict growth characteristics of successive block structures, see Figure 6 for noise in  $H_{114}(n)$. $\frac{L(n)}{n}$ has curious sign of order for $H_{114}(n)$ and $H_{146}(n)$.

\begin{figure}[!htbp]
\minipage{0.5\textwidth}
  \includegraphics[width=\linewidth]{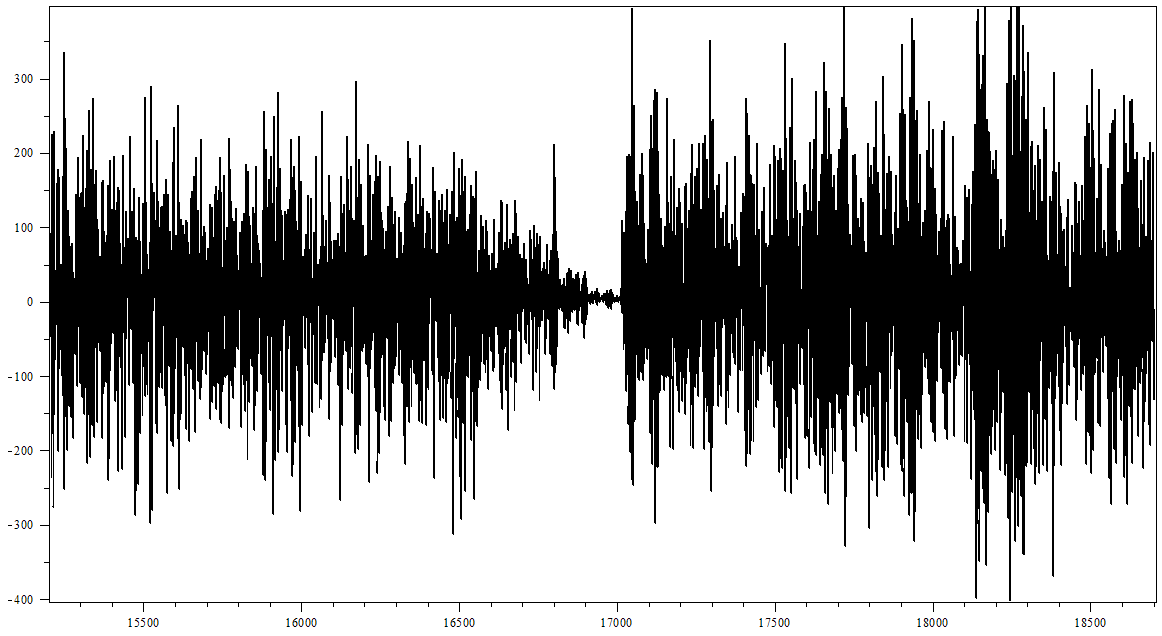}
\endminipage\hfill
\minipage{0.5\textwidth}
  \includegraphics[width=\linewidth]{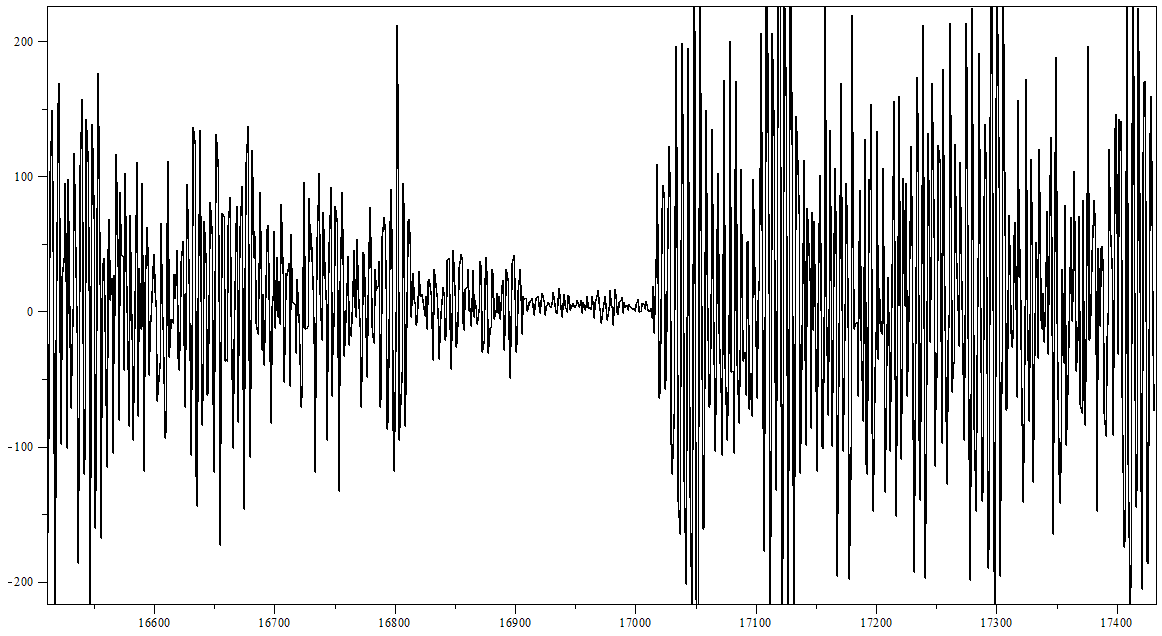}
\endminipage\hfill
\caption {Transition regions in noise of $H_{114}(n)$ for initial erratic generations.}
\label{fig:Q}
\end{figure}
\newpage
Let $S_{H_{k}}(n) = H_{k}(n)-\frac{n}{2}$ and $\big \langle S_{H_{k}}(n) \big \rangle_{t}$ denote the average value of  $S_{H_{k}}(n)$ over the $t^{th}$ generation, as described thanks to auxilary sequences for $10^{3}\leq n \leq 11\cdot10^{6}$. We can then define a variance function based on the $t^{th}$ generation:
\[
M_{t}(S_{H_{k}}(n))^2 = \big \langle S_{H_{k}}(n)^2 \big \rangle_{t} - \big \langle S_{H_{k}}(n) \big \rangle_{t}^2.
\] %, and define $\alpha(k, S_{c}(n))$ as below. 
The main aim is to compute how the standard deviations of such block structures change from one generation to the next.  We define
\[
\alpha(t, S_{H_{k}}(n)) = \log_2\!\left(\frac{M_{t}(S_{H_{k}}(n))}{M_{t-1}(S_{H_{k}}(n))}\right).
\]
%\\

\begin{table*}[!h]

\begin{center}
%\begin{ruledtabular}
\begin{adjustbox}{max width=0.85\textwidth}
\begin{tabular}{ccc}
\noalign{\smallskip}
$t$ & $\alpha(f(t), S_{H_{114}}(n))$ & $\alpha(g(t), S_{H_{146}}(n))$ \\
\noalign{\smallskip}\hline\noalign{\smallskip}

%5	&0.8402	  \\
%6	&0.7278	  \\
%7	&0.7477	  \\
1	&1.2328	&1.2712   \\
2	&1.0777	&1.1233  \\
3	&1.0534	&1.2428 \\
4	&1.1581	&1.2145  \\
5	&1.1944    &1.1708\\
6	&1.1706    &1.1908\\
7	&1.1938	&1.1988  \\
8	&1.1896	&1.1853  \\
9   &1.1955       &1.1913\\
10   &1.1910     &1.1898 \\
11   &1.1944     &1.1914\\
12   &1.1943     &1.1936\\
\hline
\noalign{\smallskip}
\end{tabular}
\end{adjustbox}
\caption {Values of $\alpha(h(t), S_{H_{k}}(n))$ for $1 \leq t \leq 12$.}
\label{tab:stat3}
%\end{ruledtabular}
\end{center}
\end{table*}

Computational results in Table 2 indicate that the $\alpha$ values oscillate around $1.19$ which suggests the mortality of $H_{114}(n)$ and $H_{146}(n)$. This resembles with $V_{c}(n)$\footnote{ $V_c$ is generated by the $V$-recurrence with the initial conditions (3,4,5,4,5,6) and last term is $V_c(3080193026) = 3101399868$.} sequence which has $\alpha$ value appears to approach $1.12$~\cite{14}. However for the $Q$-sequence and its relatives, the values of $\alpha$ approach a value around $0.88$~\cite{4,5,12,15}.

%%%%%%%%%%%%%%%%%%%%%%%%%%%%%%%%%%%%%%%%%%%%%%%%%%%%%%%%%%%%%%%%%%
\newpage
\section{Conclusion}\label{sec:conc}

In this extended abstract, an investigation of $H$-recurrence according to the initial condition patterns which are proposed by asymptotic properties of certain sequences is introduced and computational results are reported in terms of the clarification of finite chaotic solutions which this study focuses on. With the analysis of $H_{114}(n)$ and $H_{146}(n)$, growth characteristics of succesive generations and $\alpha$ values seem to have a complex relation with $r$ and $s$ values based on existence of a generational structure of chaotic solution sequence for Hofstadter-Huber recurrence relation $Q_{r,s}(n)=Q_{r,s}(n-Q_{r,s}(n-r)) + Q_{r,s}(n-Q_{r,s}(n-s))$. In particular, constructing and controlling the long chaotic sequences which have highly similar underlying regularities for $H(n)$ recurrences suggest that finding new initial condition patterns to corresponding meta-Fibonacci recurrences has potential power to new communications in discrete dynamical systems.

\section*{Acknowledgements}

Altug Alkan would like to thank Remy Sigrist and Giovanni Resta for his valuable computational assistance regarding OEIS contributions.

\section*{Conflict of Interest}

The authors declare that there is no conflict of interest regarding the publication of this extended abstract.

{\footnotesize

\renewcommand{\refname}
}


{\hfil References}
\begin{thebibliography}{99}
\bibitem{1} Tanny S.: An Invitation to Nested Recurrence Relations. Talk given at the 4th biennial Canadian Discrete and Algorithmic Mathematics Conference (Canadam). https://canadam.math.ca/2013/program/slides/Tanny.Steve.pdf (2013)
\bibitem{2} Sloane, N.J.A.: OEIS Foundation Inc.. The On-Line Encyclopedia of Integer Sequences (2019)
\bibitem{3} Hofstadter, D. \textit{Godel, Escher, Bach: an Eternal Golden Braid.} Penguin Books. (1979)
\bibitem{4} Pinn, K.: A Chaotic Cousin of Conways Recursive Sequence. Experimental Mathematics. 9, 55-66 (2000)
\bibitem{5} Pinn, K.: Order and chaos in Hofstadters Q(n) sequence. Complexity. 4, 41-46 (1999)
\bibitem{6} Ruskey, F.: Fibonacci meets Hofstadter. Fibonacci Quart. 49, 227–230 (2011)
\bibitem{7} Fox N.: Discovering Linear-Recurrent Solutions to Hofstadter-Like Recurrences Using Symbolic Computation (2017)
\bibitem{8} Fox, N.: An Exploration Of Nested Recurrences Using Experimental Mathematics. Ph.D. Thesis, Department of Mathematics, Rutgers University (2017)
\bibitem{9} Fox N.: Quasipolynomial Solutions to the Hofstadter Q-Recurrence. Integers. 16 (2016).
\bibitem{10} Balamohan B., Kuznetsov A., Tanny S.: On the behaviour of a variant of Hofstadter’s Q-sequence. J. Integer Sequences. 10, Article 07.7.1. (2007)
\bibitem{11} Dalton B, Rahman M, Tanny S. Spot-based generations for meta-Fibonacci sequences. Experimental Mathematics 20(2), 129-37. (2011)
\bibitem{12} Alkan, A., Fox, N., Aybar O.O.:On Hofstadter Heart Sequences. Complexity. 1-8 (2017)
\bibitem{13} A.S. Fraenkel, Iterated floor function, algebraic numbers, discrete chaos, Beatty subsequences,
semigroups. Trans. Am. Math. Soc. 341, 639–664 (1994)
\bibitem{14} Alkan, A., Fox, N., Aybar O.O., Akdeniz Z.: An exploration of solutions to two related Hofstadter-Huber recurrence relations, Chaos, Solitons and Fractals 138 (2020), 109900. https://doi.org/10.1016/j.chaos.2020.109900
\bibitem{15} Alkan, A.: On a Generalization of Hofstadters Q-Sequence: A Family of Chaotic Generational Structures. Complexity. 1-8 (2018) 
\bibitem{16} Alkan, A.: On a conjecture about generalized Q-recurrence. Open Mathematics. 16, 1490-1500 (2018)
\bibitem{17} Alkan A., Aybar O.O: On Families of Solutions for Meta-Fibonacci Recursions Related to Hofstadter-Conway $10000$ Sequence. In: Stavrinides S., Ozer M. (eds) Chaos and Complex Systems. Springer Proceedings in Complexity. Springer, Cham (2020)
\end{thebibliography}
\end{document}